\documentclass[11pt,a4paper]{amsart}

\newcommand{\1}{{1\!\!1}}
\usepackage{latexsym}

\usepackage{amssymb}
\usepackage{amsthm}
\usepackage{amsmath}
\usepackage{amstext}
\usepackage{epic}
\usepackage{eepic}

\numberwithin{equation}{section}

\theoremstyle{plain}
\newtheorem{thm}{Theorem}[section] 
\newtheorem{prop}[thm]{Proposition}
\newtheorem{cor}[thm]{Corollary}
\newtheorem{lem}[thm]{Lemma}

\newtheorem{theorem*}{Theorem}[]

\theoremstyle{definition}
\newtheorem{defn}[thm]{Definition}

\theoremstyle{remark}
\newtheorem{rem}[thm]{Remark}

\newcommand{\R}{\mathbb{R}}

\newcommand{\Z}{\mathbb{Z}}
\newcommand{\proj}{\mathbb{RP}}

\newcommand{\ztwo}{\mathbb{Z}_2}
\newcommand{\Hcl}{{H}^{BM}}

\newcommand{\CC}{\mathcal{C}}
\newcommand{\cc}[1]{\overline{#1}^{\mathcal{C}}}
\newcommand{\AC}{\mathcal{AC}}
\newcommand{\AS}{\mathcal{AS}}
\renewcommand{\1}{\mathbf{1}}
\newcommand{\zar}[1]{\overline{#1}^Z}
\newcommand{\G}{\mathcal{G}}

\newcommand{\la}{\longrightarrow}
\newcommand{\inv}{^{-1}}

\DeclareMathOperator{\id}{id}
\DeclareMathOperator{\sgn}{sign\,}

\DeclareMathOperator{\dist}{dist\,}

\def\accentclass@{7}
\def\makeacc@#1#2{\def#1{\mathaccent"\accentclass@#2 }}
\makeacc@\cir{017}

\title{Topology of injective endomorphisms of real algebraic sets}

\author{ Adam Parusi\'nski}

\address {D\'epartement de Math\'ematiques, U.M.R. 6093 du C.N.R.S,
Universit\'e d'Angers, 2, bd Lavoisier, 49045 Angers Cedex, France}

\email{parus@tonton.univ-angers.fr}


\begin{document}
 
\begin{abstract}
Using only basic topological properties of real algebraic sets and
regular morphisms we show 
that any injective regular self-mapping of a real algebraic set  
is surjective.  Then we show that injective morphisms between germs 
of real algebraic sets define a partial order on the equivalence classes 
of these germs divided by continuous semi-algebraic homeomorphisms.   
We use this observation to deduce that  any injective regular
self-mapping of a real algebraic set  
is a homeomorphism.  We show also a similar local property.  
All our results can be extended to 
 arc-symmetric semi-algebraic sets and injective continuous arc-symmetric
morphisms, and some results to Euler semi-algebraic sets
and injective continuous semi-algebraic morphisms.  
\end{abstract}
\maketitle


\vspace{ 1 truecm}  
In 1960 Newman \cite{newman} showed that any injective 
real polynomial map $\R^2\to \R^2$ is surjective.  This
was extended 
to the real polynomial maps $\R^n \to \R^n$ in 1962 by 
Bia{\l}ynicki-Birula and Rosenlicht \cite{birula}.  
In 1969 Ax \cite {ax} showed that any injective regular 
self-mapping of a complex algebraic variety is surjective.  
Unilike the proof of \cite{birula}, that was topological,
the proof of Ax is based on the Lefschetz principle and a
reduction to the finite field 
case.  In \cite{borel} (1969) Borel extended the idea of 
\cite{birula} and gave a topological proof of Ax' Theorem
that can be applied to 
real algebraic non-singular varieties.   The first proof of 
an analogous result 
for  singular real algebraic varieties was given in 1999
by Kurdyka \cite{kurdyka2}.  Kurdyka's proof  is based on
Borel's argument and the geometry of semi-algebraic 
arc-symmetric sets.  For more on history of the problem of 
surjectivity of injective mappings, the motivation, and 
a wide spectrum of possible applications we refer the reader 
to a recent paper of Gromov \cite{gromov}.   We note that Borel's 
proof gives as well that an injective regular 
self-mapping of a complex algebraic variety or of a non-singular
real algebraic variety is a homeomorphism.  The analogous statement 
in the general real algebraic case has not been proven till now 
(to the author's best knowledge) and doesn't not follow 
from \cite{kurdyka2}.  

In this paper we first present a new topological 
proof of Kurdyka's theorem and then extend the argument 
 to show that injective regular self-mappings of real 
algebraic varieties are homeomorphisms.  Moreover we obtain 
local versions of this result and establish a hierarchy of 
real algebraic singularities (more generally of Euler 
semi-algebraic singularities) with respect to injective 
regular mappings.  

Our approach is based on two classical topological properties of 
real algebraic sets and maps.  The first one is  Sullivan's Theorem 
that each real algebraic set is Euler (see definition \ref{Euler} below).  
The second one says that for a regular 
map $f:X\to Y$ of real algebraic sets there exists 
a proper real algebraic subset $Y'\subset Y$ such that 
the Euler characteristic of the fibres of $f$, taken 
modulo $2$, is constant on $Y\setminus Y'$.  This is the crucial 
observation for the problem since it implies that the image of injective 
regular map of real algebraic varieties is algebraically
constructible, i. e. belongs to the Boolean algebra 
generated by real algebraic sets.  

In section \ref{TP} we recall 
basic properties of Euler semi-algebraic sets.  In particular, 
these sets have a well-defined fundamental class in the
Borel-Moore homology with $\Z_2$ coefficients.  As an easy corollary, 
proposition \ref{compact}, we obtain that any injective 
continuous semi-algebraic self-map of 
a compact  Euler semi-algebraic set is a homeomophism and an 
analogous local version.  The class of algebraically constructible
sets is introduced in section \ref{AC}.  Its topological properties, 
similar to the ones of real algebraic sets are established in sections 
\ref{AC} and  \ref{Examples}.  The sujectivity theorem of Kurdyka is
proven in section \ref{Borel-Kurdyka}.  The proof is based on Borel's
argument in non-singular case and the hunt for invariant subsets.  More
precisely, let $X$ be a real algebraic set and let $f:X\to X$ be
injective regular.  We show that there is a subset 
$\tilde \Sigma \subset X$, of dimension smaller than $X$, such that 
$f(\tilde \Sigma) \subset \tilde \Sigma$ and $U = X\setminus  \tilde
\Sigma$ is a topological manifold.  We do not know whether such
$\tilde \Sigma$ can be found algebraic, but the set we construct is
algebraically constructible and closed in $X$.  Then we apply the
induction on dimension to $f|_{\tilde \Sigma} : \tilde \Sigma\to
\tilde \Sigma$ to show that it is surjective (this step shows that  we
have to extend the category we work with from the algebraic sets to 
the algebraically constructible ones).  Then $f(U) \subset U$ and
we conclude that $f(U) = U$ by Borel's argument.  

In section
\ref{Hierarchy} we study germs of 
real algebraic sets (or more generally of Euler semi-algebraic sets)
and injective mappings between them.  For two such germs we write 
$(Y,y) \prec (X,x)$ if there is a germ $f:(Y,y) \to (X,x)$, $f$ being  
continuous, injective, and semi-algebraic.  We show  that
$\prec$ is a partial order (on the set of germs up to 
semi-algebraic homeomorphisms).  One may think that  $(Y,y) \prec
(X,x)$ means that the singularity type of $(Y,y)$ is less complicated
than that of $(X,x)$.  There are two precise result that justify this
interpretation.  Firstly a non-singular germ is minimal (but not 
necessarily unique as such) among the germs of
a fixed dimension, corollary \ref{invariance2}.  
Secondly, if $x\in X$ then  the singularity type of $(X,x')$ for   
$x'$ close to $x$  cannot be more complicated than that 
of $(X,x)$, i. e. if $(X, x) \prec (X,x')$ then both germs are
homeomorphic.  Using this hierarchy of singularities we show the main
result that any
injective (and so bijective) regular self-map of real algebraic set 
is a homeomorphism, theorem \ref{homeo2}.  
(Surprisingly, in the case of sets of
pure-dimension, the assumption can be weaken considerably: any
bijective continuous semi-algebraic self-map of a pure-dimensional 
Euler semi-algebraic set is a homeomorphism, cf. theorem
\ref{homeo1}.)

In the paper we formalize what properties of algebraically
constructible sets we need by introducing the notion of 
constructible category of semi-algebraic sets, definition
\ref{ACdefinition}, and 
then we proceed in any such category.    
Though this makes some proofs, e. g.  of theorem \ref{ACC}, more
involved it has the advantange of extending automatically the results 
of the paper to semi-algebraic arc-symmetric sets, cf. 
subsection \ref{AS} .  It gives
also a bridge to the approach of Kurdyka, 
though we use only basic topological properties of arc-symmetric 
sets.  Algebraically constructible sets form the smallest
constructible category, arc-symmetric sets the biggest one, see
section \ref{Remarks}.  In section \ref{Examples} we relate these categories
to the integration with respect to the Euler characteristic and the 
algebraically constructible and the Nash-constructible functions.  

\smallskip \emph{Notation and terminology.}  In this paper a ``regular
set''
means a pure-dimensional topological manifold and a ``stratification''
means a finite decomposition into such regular, not necessarily
connected, sets.  The Zariski closure of $X$ will be denoted by $\zar
X$.

\bigskip
\section{Topological preliminaries}
\label{TP}
\medskip

By \emph{a semi-algebraic set} we mean a semi-algebraic subset of an 
affine space $\R^N$, or of a projective space $\proj^N$.  We often,
but not always, assume that a semi-algebraic set is locally 
compact, that is locally closed as a subset of $\R^N$ or 
$\proj^N$.  
A \emph{semi-algebraic map} is a 
map between two semi-algebraic sets with semi-algebraic
graph.  Semi-algebraic maps may be continuous or not.  
Similarly, by \emph{algebraic sets} we mean the algebraic subsets
of affine or projective spaces. Our discussion can be extended to 
semi-algebraic and algebraic subsets of real algebraic varieties
that we will not do for simplicity of exposition.  
For the background on semi-algebraic and real algebraic
sets we refer the reader to \cite{BCR}.  

We say that a semi-algebraic set $X$ is \emph{(topologically) 
regular of dimension $d$}  
at $x\in X$ if there is a neighbourhood of $X$ in $x$ that is a 
topological manifold of dimension $d$. In this paper we shall 
drop the word ``topologically'' and call 
topologically regular sets regular.   We say that $X$ is
\emph{regular at $x$} if it is regular of dimension 
$d=\dim X$. In this case we also say that $x$ is a
\emph{regular point of $X$}.  
$X$ is \emph{regular} if it is regular at all its
points.

In this section we discuss basic topological properties of Euler
semi-algebraic sets.  Let $X$ be semi-algebraic and $x\in X$.  
The local Euler-Poincar\'e characteristic of $X$ at $x$ is, 
by definition, 
$$
\chi(X,X\setminus x) = \sum_i (-1)^i \dim H_i(X, X\setminus x; 
\ztwo )  . 
$$

\begin{defn}\label{Euler}
We call a locally compact semi-algebraic set $X$,  
\emph{Euler} 
if $\chi(X,X\setminus x)$ is odd for all $x\in X$. 
\end{defn}
 
Let $X\subset \R^N$ and let $S(x, \varepsilon)$ be 
the sphere
of radius $\varepsilon>0$ centered at $x$.  By the 
local conic structure lemma, cf. \cite{BCR} (9.3.6), for 
$\varepsilon$ sufficiently small the topological type of 
$S(x,\varepsilon)\cap X$ is independent of 
$\varepsilon$.  This space is called the \emph{link} of $x$ in
$X$.  It is easy to check that $\chi(X,X\setminus x)$ is odd 
iff the Euler characteristic of the link of $X$ at $x$ is even.  
By Sullivan's theorem \cite{sullivan} each real algebraic set 
is Euler.  

Let $X$ be a locally compact semi-algebraic set.  
Consider the Borel-Moore homology $\Hcl (X)$ of $X$,  
see \cite{borelmoore}.  
Since $X$ can be triangulated the
Borel-Moore homology of $X$ is isomorphic to the 
simplicial homology of
$X$ with closed supports (\emph {i.e.} using possibly 
infinite simplicial chains). If 
$X'$ is a compact space containing $X$ such that 
$X'$ is triangulated with
$X'\setminus X$ as a subcomplex, then the Borel-Moore 
homology of $X$ is isomorphic to the relative simplicial 
homology of the pair
$(X', X' \setminus X)$.  In particular if 
$\bar X = X\cup \{\infty\}$  is 
 the one point compactification of $X$, 
that can be identified with a semi-algebraic set, cf.  
\cite{BCR} (2.5.9), we have  
$$
\Hcl _i (X; \ztwo ) = H_i (\bar X, \{\infty\}; \ztwo ). 
$$ 
For an open semi-algebraic $U\subset X$  there is 
a restriction homomorphism 
\begin{equation}\label{open}
 \Hcl _i (X; \ztwo ) \to \Hcl_i (U;\ztwo ). 
\end{equation}
Then $Y=X\setminus U$ is a closed semi-algebraic subset of $X$ 
and the exact sequence of a triple for classical homology 
gives an exact sequence, 
cf. \cite {BCR} (11.7.12),
\begin{equation}\label{exactBM}
\cdots \to \Hcl_{i+1} (U; \ztwo) \to \Hcl_i (Y; \ztwo) \to 
\Hcl_i (X; \ztwo) \to \Hcl_{i} (U; \ztwo) \to \cdots .  
\end{equation}

Semi-algebraic sets can be triangulated, a compact 
semi-algebraic set $X$ is homeomorphic to a finite 
simplicial complex.  Let $\dim X = n$.  If $X$ is Euler then for 
each $(n-1)$-dimensional simplex $\tau$ 
there is an even number of $n$-dimensional 
simplices that contain $\tau$ (we shall consider all simplices closed).  
This shows that the formal sum of $n$-simplices forms a 
$\ztwo$-cycle.    Its class in ${\rm H}_n(X;\ztwo )$ is non-zero.  
It is  called
\emph{the fundamental class of $X$} and denoted by $[X]$. It is 
characterized by the property that its image under 
${\rm H}_n(X;\ztwo )\to {\rm H}_n(X, X\setminus x;\ztwo )$ is 
nonzero for all $x$ regular in $X$.  The fundamental 
class is independent of 
the choice of triangulation.   

If $X$ is locally compact and Euler 
then  $\bar X$ is Euler.  
This is a simple exercise.  We provide a proof for 
the reader convenience.

\begin{prop}\label {compactification}
The one point compactification of an Euler semi-algebraic 
set is Euler.
\end{prop}

\begin{proof}
Triangulate $\bar X$ so that $\{\infty \}$ is a simplex.  
A triangulated space is Euler if for each simplex $\sigma$ 
the number of simplices $\tau$ such that 
$\sigma \subset \tau, \sigma \neq \tau$,  is even.  
Note that each k-simplex  $\tau$ contains exactly $2^{k+1} -2$ 
non-empty simplices $\sigma$ such that $\sigma \neq \tau$.  
 Thus, if we count the number 
of relations $\sigma \subset \tau, \sigma \neq \tau$,  among 
non-empty simplices,  we get an even number.  
Hence there is always an even number of simplices $\sigma$ 
having an odd number of adherent simplices $\tau$, i. e. 
$\sigma \subset \tau, \sigma \neq \tau$.  
This shows the proposition since for the given 
triangulation of $\bar X$ all simplices except maybe 
$\{\infty\}$ have an even number of adherent simplices 
by assumption.  
\end{proof}

For a locally compact semi-algebraic set of dimension $n$ 
we define the \emph{fundamental class } of 
$X$ as the 
image of $[\bar X]$ by the canonical homomorphism, 
cf. \cite{BCR} (11.4.1), 
$$
H_n(\bar X; \ztwo ) \to 
H_n (\bar X, \{\infty \}; \ztwo ) = \Hcl _n (X; \ztwo ). 
$$ 

Here are simple general observations related to 
the problem of surjectivity of injective mappings. 

\begin{prop} \label{compact}
\item [(i)]
Let $f:X\to X$ be a continuous semi-algebraic self-map of 
an Euler compact semi-algebraic $X$.  If $f$ is injective 
then it is a homeomorphism.  
\item [(ii)]
Let $f:(X,x)\to (X,x)$ be the germ of a continuous 
semi-algebraic self-map of an Euler semi-algebraic $X$.  
Then if $f$ is injective then it is the germ of a local
homeomorphism.  
\end{prop} 

\begin{proof}
To prove (i) it suffices to show that $f$ is surjective.  
Suppose that this is not the case.  The image $Y= f(X)$ of $f$ 
is a closed semialgebraic subset of $X$ and 
$f: X \to f(X)$ is a homeomorphism.  Let 
$U = X\setminus Y$ and let $d= \dim U$.  By \eqref{exactBM} 
the following sequence is exact 
\begin{equation}\label{}
0 \to  \Hcl_d (Y; \ztwo) \to 
\Hcl_d (X; \ztwo) \to \Hcl_{d} (U; \ztwo) \to \cdots .  
\end{equation}
But $X$ and $Y$ are homeomorphic and hence 
$\dim \Hcl_d(X;\ztwo ) = \dim \Hcl _d(X;\ztwo )$.  
This would imply that 
$\Hcl_d (X; \ztwo) \to \Hcl_{d} (U; \ztwo)$ is the zero 
homomorphism that is impossible by the following lemma.   

\begin{lem}\label{openrestriction}
Let $X$ be an Euler locally compact semi-algebraic set 
and let $U\subset X $ be open
semi-algebraic.  Then the fundamental class of $U$ is 
in the image of restriction homomorphism   
$\Hcl _{\dim U} (X; \ztwo ) \to \Hcl_{\dim U} (U;\ztwo )$. 
\end{lem}

\begin{proof}
If  $\dim U = \dim X$,  then it is easy to see that  
$[U]$ is the image of $[X]$.  In the general case let $d = \dim U$.  
Let $\bar X$ be the one point compactification of $X$.  
Triangulate 
$\bar X$ so that $\bar X\setminus U$ is a subcomplex and 
consider the first 
barycentric subdivision of $\bar X$.  Then, by \cite{sullivan}, see 
also \cite{akin} proposition 1, 
the formal sum of d-simplices is a $\ztwo$-cycle.  
(It represents the $d$-th
homology Stiefel-Whitney class of $\bar X$.)  Its image in 
$H_d(X, X\setminus U; \ztwo) = \Hcl _d(U;\ztwo )$ is the fundamental 
class of $U$.  This shows the lemma.    
\end{proof}

Let $B(x, \varepsilon) $ be an open ball in $X$ 
centered at $x$ and of radius $\varepsilon>0$.  By the local conic
structure,  \cite{BCR} (9.3.6), 
$$
\Hcl _i (B(x, \varepsilon) ;\ztwo ) 
\simeq H_i(X, X\setminus x;\ztwo ),    
$$
for $\varepsilon>0$ and small.  
This holds not only for small balls but for any semi-algebraic family 
of open neighbourhoods of $x$ in $X$ in particular for the familly 
$f\inv (B(x, \varepsilon))$, $\varepsilon >0$. That is we claim that 
\begin{equation}\label{local2}
\Hcl _i (B(x, \varepsilon) ;\ztwo ) = 
\Hcl _i (f\inv (B(x, \varepsilon)) ;\ztwo ) ,    
\end{equation}
for $\varepsilon>0$ and small. 
Indeed, this is can be obtain by the following classical argument.  Let 
$\rho: X\to [0,\infty)$ be given by  
$\rho (y) = \dist (f(y), x)$.  Denote $U(x, \varepsilon) = \rho \inv
([0,\varepsilon))$. Then, by topological triviality of $\rho$ over 
$(0,\varepsilon_0)$, for an $\varepsilon_0>0$, $U(x, \varepsilon)
\subset U(x, \varepsilon')$ is a homotopy equivalence if 
$0<\varepsilon \le \varepsilon' \le \varepsilon_0$.  Considering the 
homomorphism of cohomology groups induced by the 
inclusions 
$$
B(x,\varepsilon_1) \subset  U(x, \varepsilon_2) \subset 
B(x,\varepsilon_3)\subset  U(x, \varepsilon_4), \quad 
0< \varepsilon_1 \ll \varepsilon_2 \ll  
\varepsilon_3 \ll \varepsilon_4 \ll  \varepsilon_0, 
$$
we get $\Hcl _i (B(x, \varepsilon_1) ;\ztwo ) = 
\Hcl _i ( U(x, \varepsilon_2) ;\ztwo )$ and hence \eqref{local2}.  
(Alternatively one may show $\Hcl _i (U(x, \varepsilon)) ;\ztwo )
\simeq H_i(X, X\setminus x;\ztwo )$ using the local conic structure 
induced of $X$ at $x$ induced by $\rho$.)

Now  (ii) follows from the same argument as (i) applied to $B(x,
\varepsilon)$, its closed subset $f(X) \cap B(x,
 \varepsilon)$, and $U= B(x, \varepsilon) \setminus f(X)$.  
\end{proof}

\begin{rem}\label{codimension1}
Note that for the existence of fundamental class it 
is not necessary to assume tha $X$ is Euler.  It suffices that 
$X$ is \emph{Euler in codimension one} that is there 
 exists a semi-algebraic subset $Y\subset X$, $\dim Y \le 
\dim X-2$, such that $\chi(X,X\setminus x)$ is odd for every 
$x\in X\setminus Y$.  (This condition is actually equivalent to the 
existence of fundamental class.)  
Thus, in particular, proposition \ref{compact} and lemma 
\ref{openrestriction} hold for 
$X$ Euler in codimension 1 provided it is of pure 
dimension (then the set $U$ of the proof has the same 
dimension as $X$).  
\end{rem}

We shall need later the following simple observation.  

\begin{lem}\label{invariance} 
 Let $X$ be a connected topological manifold and a semi-algebraic 
set and let 
$Y$ be a closed semi-algebraic subset of $X$.  If $Y$ is 
Euler in codimension 1 then either $Y=X$ or $\dim Y< \dim X$.  
\end{lem} 

\begin{proof}
Let $n=\dim X$.  
Suppose $U= X\setminus Y$ is non-empty.  Then 
\eqref{exactBM} gives an exact sequence 
\begin{equation}
0 \to \Hcl_n (Y; \ztwo) \to 
\Hcl_n (X; \ztwo) \to \Hcl_{n} (U; \ztwo) \to \cdots .  
\end{equation}
By assumption the only non-zero element of $\Hcl_n (X;\ztwo )$ is 
$[X]$.  Its image  in $\Hcl_n (U;\ztwo )$  equals 
$[U]$ by lemma \ref{openrestriction} so is non-zero and hence   
$\Hcl_n (Y; \ztwo) =0$.   But if $Y$ is Euler in codimension 1 then 
$\Hcl_{\dim Y} (Y; \ztwo) \neq 0$ and therefore $\dim Y < n$. 
\end{proof}

\begin{cor}\label{invariance2} 
 Let $f:(Y,y)\to (X,x)$ be an injective continuous and semi-algebraic
 map, where $(X,x), (Y,y)$ are germs of locally compact semi-algebraic
 sets, $\dim_x X = \dim _y Y$.
 Suppose that  $(X,x)$ be a topological manifold and that $(Y,y)$ 
is Euler in
 codimension 1.  Then $f$ is a local homeomorphism.    
\end{cor}

\medskip
\section{Constructible Categories}
\label{AC}
\medskip

First recall two basic topological properties of 
real algebraic sets and regular maps.    

\begin{thm} \label{sullivan}{\rm (Sullivan, \cite {sullivan}) }
Every real algebraic set is Euler.  
\end{thm}

Regular maps satisfy the following basic
topological mod $2$ property (cf. \cite{akbking} (2.3.2)).  

\begin{thm} \label{standard} {\rm (see for instance 
\cite{akbking} (2.3.2))} 
Let $F:X\to Y$ be a regular map of real algebraic sets.
Then there is a proper algebraic subset $Y'\subset Y$ such 
that the Euler characteristic $\chi (F\inv (y)) $, $y\in Y$, 
is constant modulo 2 on $Y\setminus Y'$.    
\end{thm} 

\begin{defn} \label{ACdefinition}
We say that  a semi-algebraic set $X$ of $\R^N$  or 
$\proj^n$ is \emph{algebraically constructible} 
(or Zariski constructible) if $X$ belongs to the Boolean algebra
generated by the algebraic subsets of $\R^N$.  
Then we write $X\in \AC$ for short. 
\end{defn} 

The algebraically constructible sets satisfy similar topological 
properties to the real algebraic sets.  For instance it can 
be easily showed by Sullivan's Theorem that each locally 
compact algebraically constructible set is Euler.  But they 
form a strictly wider class than the algebraic sets.  
Closed or even compact algebraically constructible sets 
are not necessarily algebraic.  Moreover the image of a regular
injective algebraic morphism is not necessarily algebraic but it is 
algebraically constructible.  Indeed, let $F:X\to Y$  be  
a regular injective map of real algebraic sets.  Then the 
Euler characteristic of fibers of $F$ is either $0$ or $1$.  
Thus, if $Y'\subset Y$ is the subset given by theorem 
\ref{standard}, then $Y\setminus Y'$ is either contained 
in or disjoint with $F(X)$.  Then, by 
inductive argument on $\dim Y$ and the number of irreducible
components of $Y$, applied to the restriction 
$F|_{F\inv (Y')} : F\inv (Y') \to Y'$, 
we obtain that $F(X)$ is algebraically constructible.

\medskip
\begin{defn}
Let $\CC$ be a sub-collection of semi-algebraic sets.  We call
a map between two semi-algebraic sets of $\CC$ 
a \emph{$\CC$-map} if its graph is in $\CC$.   
We say that $\CC$ is a \emph{constructible category} 
if it satisfies the following axioms:

\begin{enumerate}
\item [A1.] 
$\CC$ contains algebraic sets.  
\item [A2.]
$\CC$ is stable by boolean operations
$\cap, \cup, \setminus$.
\item [A3.]
\begin{enumerate}
\item 
$\CC$ is stable by the inverse images by $\CC$-maps.  
\item 
$\CC$ is stable by images by injective $\CC$-maps.  
\end{enumerate}
\item [A4.]
For each locally compact $X$ in $\CC$ there exists 
a semi-algebraic subset $Y\subset X$, $\dim Y \le 
\dim X -2$, such that $X\setminus Y$ is Euler.  
\end{enumerate}
\end{defn}

Recall after section \ref{TP} that the last axioms means 
that $X$ has a well-defined
fundamental class in the Borel-Moore homology with
$\Z_2$ coefficients.  

If a semi-algebraic $X$ belongs to a constructible category
$\CC$ then we say that $X$ is \emph {$\CC$-constructible}
and write $X\in \CC$.  
We shall show in the Section \ref{Examples} below that
the algebraically constructible sets form a constructible
category.  This is the smallest constructible category.
The arc-symmetric sets,
see Section \ref{Examples} below, form the largest constructible 
category.  Using this observation we shall show,  
see subsection \ref{Remarks}, that the axioms A1-A4 
actually imply that each 
locally compact set in a constructible category is 
Euler.  Constructible categories, as the 
 algebraically constructible sets or the 
arc-symmetric sets, share many other interesting 
properties, as for instance the existence of 
$\CC$-closure.  

\begin{thm}\label{closure}
Let $\CC$ be a constructible category.  Then the  
semi-algebraic sets posses a well-defined
closure in $\CC$.  That is for any 
given locally compact $\CC$-constructible set $X$
and a semi-algebraic
$A\subset X$ there is the smallest closed subset of $X$
that belongs to $\CC$ and contains $A$.  We denote it
by $\cc A$.  
Any other closed subset of $X$ that is in $\CC$ and
contains $A$ must contain $\cc A$. 
\end{thm}

\begin{proof}
Induction on $\dim X$.  

Consider first the case $\dim (X\setminus A) < \dim X$.
Replacing $A$ by its topological closure we may suppose $A$ closed in $X$.  
Let $S= \zar {(X\setminus A)}$ denote the Zariski closure of 
$X\setminus A$.  Then $X' = S\cap X$
is $\CC$-constructible.  Set $A' = S\cap A$ and define
$$
\cc A = A \cup \cc {A'} = (X\setminus S) \cup \cc{A'},
$$
where $\cc {A'}$ is the $\CC$-closure of $A'$ in $X'$ that
exists by the inductive hypothesis.  The first union shows that 
$\cc A$ is closed and the second that it is $\CC$-constructible. 
 One may easily check
that $\cc A$ satisfies the required minimality 
property.

Define $Z:= \zar A$, 
where  $\zar A$ denote the Zariski closure of $A$.
Replacing $X$ by $X\cap Z$, if necessary, we may suppose
$A\subset X\subset Z$ and $Z = \zar X$.  There
is an algebraic subset $S\subset Z$ such that
$\dim S < \dim Z$, $Z\setminus S$ is a topological 
manifold of pure dimension,
and $X\setminus S$ is the union of some connected
components of $Z\setminus S$.
To each $C\subset X$, $C$ closed in $X$,
$\CC$-constructible and containing
$A$, we
associate the set of connected components of
$X\setminus S$ that are entirely contained in $C$.  Taking the
intersection of finitely many such $C$ we find such
a set, denoted $C_0$, that has minimal number
of such components.  Let $A_1$ denote the union of
these components.  Thus all closed in $X$, $\CC$-constructible
sets containing $A$ contain $A_1$.  

Let $X'= (X\cap S) \cup C_0$.  $X'$ is $\CC$-constructible and 
closed in $X$.  
By lemma \ref{invariance}, 
$\dim (C_0 \setminus (A\cup A_1)) <\dim C_0= \dim X$ and
hence $\dim (X' \setminus (A\cup A_1)) <\dim X' = \dim X$.  
Then the $\CC$-closure of $A\cup A_1$ in $X'$,
that exists by the first case considered above, is the
$\CC$-closure  of $A$ in $X$.  This ends the proof.  
\end{proof}

\begin{rem}\label{remark}
Clearly $\dim \cc A = \dim A$.  If $A$ is $\CC$-constructible 
then $\cc A = A\cup \cc {(\bar A\setminus A)}$ and hence 
$\dim (\cc A \setminus A) < \dim A$. 
\end{rem}

\begin{thm}\label{ACC}
Algebraically constructible sets form a constructible 
category.
\end{thm}

This follows by standard argument from theorems 
\ref{sullivan} and \ref{standard} but we postpone 
the proof until section \ref{Examples} where 
we deduce it from the topological properties of 
algebraically constructible functions (a formal 
approach to the integration with respect to the 
Euler characteristic on real algebraic sets).  Let us 
just note that A1 and A2 are satisfied trivially and 
A4 follows from theorem \ref {sullivan} by the additivity 
of the Euler characteristic with compact supports 
(on the link).  A3 for regular mappings is easy 
(we have showed already A3 (b)). 

In particular there is the algebraically constructible 
closure.  
Let $A$ be a semi-algebraic subset of a real
algebraic set $X$.  Then there exists the smallest set, 
denoted by $\overline A^{\AC}$,
closed  and algebraically constructible in $X$ that 
contains $A$.
Clearly $\overline A \subset \overline A^{\AC} \subset
\zar A$, where $\zar A $ denote Zariski closure of $A$, 
 but each inclusion may be strict.  
For instance if $A$ is the regular part of the Whitney umbrella
$\{zx^2=y^2\} \subset \R^3$ then $\overline A^{\AC} =
\zar A$ and is strictly bigger than $\overline A$.  If $A$ is the 
regular part of the Cartan umbrella
$\{z(x^2+y^2)=x^3\} \subset \R^3$ then
$\overline A  = \overline A^{\AC}$ but is strictly smaller
than $\zar A$.


\bigskip
\section{Proof of Borel-Kurdyka Theorem}
\label{Borel-Kurdyka}
\medskip

First we recall the theorem of Borel \cite{borel} restated
for constructible categories.  

\begin{thm}\label{borel}
Let $X$ be a locally compact
semi-algebraic set belonging to a constructible
category $\CC$.  
Let $f:X\to X$ be an injective continuous $\CC$-map.
If $f$ is open then it is a homeomorphism.
In particular such $f$ is surjective.    
\end{thm}

The original theorem of Borel is stated for non-singular real
algebraic sets and injective regular maps.  Such maps
are open by the invariance of domain, cf. \cite{munkres}  
(one may use alternatively corollary \ref{invariance2}).   
The openess is also guaranteed if $X$ is a normal complex 
algebraic variety and 
$f:X\to X$ is an injective regular morphism of complex 
algebraic varieties.  Then passing to normalization Borel 
shows that any injective regular self-mapping of a 
complex algebraic variety is surjective.    

We shall recall below the main idea of \cite{borel}.  
The existence of the fundamental class in the Borel-Moore homology 
with $\Z_2$ coefficients plays the crucial role in this argument.     
Let $d=\dim (X\setminus f(X))$.  
Denote $f^{k+1} = f\circ f^k$, $f^0 =\id_X$, and 
$Y_k = X\setminus f^k(X)$. The sets $Y_k$ are closed in 
$X$ and Euler in codimension 1, cf. remark \ref{codimension1}.  
The set  $Y_\ell\setminus Y_{\ell-1}$ are all disjoint and 
of dimension $d$ and hence $[Y_\ell]-[Y_{\ell-1}] =
[Y_\ell \setminus Y_{\ell-1}]$, $\ell = 1, \ldots , k$,  
are linearly independent in $\Hcl_d(Y_k,\Z_2)$.  Thus 
\begin{equation}\label{lowerbound}
\dim \Hcl_d (Y_k,\Z_2) \ge k~.
\end{equation}
Note that $f^k : X \to f^k(X)$ is a homeomorphism and 
$f^k(X)$ is open in $X$.  Thus the long exact sequence  
\eqref{exactBM},
$$
\la \Hcl_{d+1} \big( f^k(X),\Z_2\big) \la \Hcl _d(Y_k,\Z_2) \la
\Hcl_d (X,\Z_2) \la, $$
gives a bound independent on $k$ 
$$\dim \Hcl_d(Y_k,\Z_2) \le \dim \Hcl_d(X,\Z_2) + \dim
\Hcl_{d+1} (X,\Z_2), $$
that contradicts \eqref{lowerbound}.

\begin{thm}\label{mainthm}
Let $X$ be a locally compact
semi-algebraic set belonging to a constructible
category $\CC$.  
Let $f:X\to X$ be a continuous
morphism in $\CC$ that is injective as a map.
 Then $f$ is surjective.  
\end{thm}

By theorems \ref{ACC} and \ref{ASC}, proven in the next
sections, we get the following. 

\begin{cor}\label{borelkurdyka} (Borel-Kurdyka Theorem)
Any injective regular 
self-mapping of a real algebraic variety is surjective.
\end{cor}

\begin{cor} (cf. \cite{kurdyka2}) 
Any injective continuous self-mapping of a locally 
compact arc-symmetric set, 
the graph of which is arc-symmetric, is surjective.
\end{cor}

\begin{proof}[Proof of theorem \ref{mainthm}] 
We proceed by induction on $n= \dim X$.  
Let $\Sigma = Sing(X)$ be the set of topological singularities 
of $X$.  That is $\Sigma$ is the disjoint union of the 
following sets 
\begin{eqnarray*}
& S & = \{ x\in X ; \dim _x X= n, (X,x) \text { not homeomorphic to } 
(\R^n,0) \}, \\ 
& A & = \{ x\in X ; \dim _x X < n,  \}. 
\end{eqnarray*}    
By the existence of locally topologically trivial
semi-algebraic
stratification of $X$ both $A$ and $S$ are semi-algebraic.  It is 
easy to see that $S$ and $\Sigma = S\cup A$ are closed in $X$. 
We shall show that 
\begin{equation}\label{ccSigma}
f(\cc \Sigma) \subset \cc \Sigma.
\end{equation}

 By lemma \ref{invariance}, $f(S) \subset S$.  Therefore, 
by the lemma below, $f(\cc S) \subset \cc S$. 

\begin{lem}
Let $B\subset X$ such that $f(B)\subset B$.  Then 
$f(\cc B) \subset \cc B$. 
\end{lem}

\begin{proof} 
$f(B) \subset B$ gives
\begin{equation} \label{inclusions}
B\subset f \inv (f(B)) \subset f\inv (B) \subset 
f\inv (\cc B).
\end{equation}
Note that $f\inv (\cc B)$ is closed and in
$\CC$.  Hence \eqref{inclusions} gives
$\cc B \subset f\inv (\cc B)$.  That is 
$f(\cc B) \subset f(f\inv (\cc B)) \subset \cc B$
as required. 
\end{proof} 

Since $\dim \cc S<\dim X$ and by inductive assumption  applied to 
$f_{| \cc S}: \cc S\to \cc S$ 
\begin{equation}\label{ccS}
f(\cc S) = \cc S .
\end{equation}
Let $U = X\setminus \cc \Sigma$.  We shall show $f(U)\subset U$.  
Firstly, by \eqref{ccS} and the 
definition of $A$,  $f(U) \subset X\setminus (A\cup \cc S)$. 
Since  $A\cup \cc S$ is closed in $X$, $X\setminus (A\cup \cc S)$ is 
open in $X$.   By Lemma \ref{invariance2}, $f(U)$ is open in 
$X\setminus (A\cup \cc S)$ and hence in $X$.  Consequently 
$X\setminus f(U)$  is  closed in $X$  and, by axiom A3 (b),
$\CC$-constructible.  Thus, since   
$A\cup \cc S \subset X\setminus f(U)$, 
by the definition of $\CC$-closure, 
$\cc \Sigma = \cc A \cup \cc S \subset X\setminus f(U)$.
That is $f(U) \subset U$ as claimed.  Since $U$ is a $\CC$-constructible 
topological manifold and $f_{|U}:U \to U$ is injective, 
it is surjective by Theorem \ref{borel}. 
This means  $f(U) = U$ and  shows \eqref{ccSigma}.  
By the inductive assumption, 
$ f(\cc \Sigma) = \cc \Sigma $,  and therefore $f(X)=X$ as claimed.  
\end{proof}

We show in section \ref{Hierarchy} below that $f$ is a homeomorphism.  
This will give in particular that  $f(A)=A$, $f(S)=S$, and 
$f(\Sigma) = \Sigma$.


\bigskip
\section{Examples of Constructible Categories}
\label{Examples}
\medskip

\subsection{Algebraically Constructible Sets}

The algebraically constructible sets, see definition 
\ref{ACdefinition}, can be studied using the theory of 
algebraically constructible functions, 
see \cite{mccroryparusinski}.  
Let $Z$ be a real algebraic set. An \emph{algebraically
constructible function} on $Z$ is an integer-valued function
which can be written as a finite sum
\begin{equation}\label{acf}
\varphi = \sum m_i
f_{i*}\1 _{Z_i},
\end{equation} where for each $i$, $Z_i$ is an
algebraic set, $\1 _{Z_i}$ is the characteristic function of
$Z_i$, $f_i:Z_i\to Z$ is a proper regular map, and $m_i$ is an
integer.

\begin{thm}\label{ACequiv}
Let $X$ be a semi-algebraic subset of $\R^N$.  
The following conditions are equivalent:  
\begin{enumerate}
\item [(i)] 
 $X$ is an algebraically constructible set.
\item [(ii)]
The characteristic function $\1 _X$ of $X$ is an algebraically
constructible on $\R^N$.  
\item [(iii)] 
There exist polynomials $g_1, \ldots, g_k \in
\R[x_1,\ldots, x_N]$ such that $\1 _X = \sum \sgn g_i$.
\item [(iv)]
There exists a real algebraic morphism $F: Z\to \R^N$ such that
\begin{eqnarray*}
& x\in X \text { iff } \, \chi (F\inv (x)) = 1 \,\\
& x\notin X \text { iff } \, \chi (F\inv (x)) = 0 .
\end{eqnarray*} 
\item [(v)]
There exists a real algebraic morphism $F: Z\to \R^N$ such that
\begin{eqnarray*}
& x\in X \text { iff } \, \chi (F\inv (x)) \equiv 1 \pmod 2\,\\
& x\notin X \text { iff } \, \chi (F\inv (x)) \equiv 0 \pmod 2 .
\end{eqnarray*} 
\end{enumerate}
\end{thm}

\begin{proof}
(i) implies (ii) by the definition of algebraically
constructible functions.  (ii) is equivalent to (iii) by
\cite{parusinskiszafraniec1} or \cite{parusinskiszafraniec2}.  
(ii) (resp. (iii)) implies
(iv) by \cite{mccroryparusinski}
(resp. \cite{parusinskiszafraniec1}).  (iv) implies (v)
trivially.  
Finally (v) implies (i) follows easily from Theorem \ref{standard}. 
\end{proof}


\begin{proof}[Proof of theorem \ref{ACC}] 
We show it only for subsets of $\R^N$.  The proof for the 
subsets of $\proj^N$ is similar.  

Axioms A1 and A2 are satisfied trivially.  
To show A3 we use the fact that the 
algebraically constructible functions 
are stable  by the inverse image and the push forward by 
regular maps, see \cite{mccroryparusinski}.   
Let $f:X\to Y$ be a map between two algebraically 
constructible subsets of $\R^n$ and $\R^m$ respectively 
and suppose that the graph of $f$, 
$\Gamma_f\subset \R^n\times \R^m$, belongs to $\AC$.  
Denote the projections of $\Gamma _f$ on $X$ and $Y$ by 
$\pi_X$ and $\pi_Y$ respectively.  Let $B\subset Y$, $B\in \AC$.  
Then $\pi_Y\inv (B)\in \AC$ and hence the characteristic 
function of $\pi_Y\inv (B)$ is algebraically constructible.  
Consequently its push-forward onto $X$ by $\pi_X$,  
that equals  the characteristic function
of $f\inv (B)$, is algebraically constructible.  
This shows that $f\inv (B) \in \AC$ and A3 (a) is shown.  
If $A\subset X$, $A\in \AC$, and $f|_A$ is 
injective then $\pi_Y$ restricted to $\pi_X\inv (A)$ is injective 
and the same argument as above shows that 
$f(A) = \pi_Y(\pi_X\inv (A))\in \AC$.  This ends the proof of 
A3.  
  
A4 follows from the fact  the link of an algebraically 
constructible function is even \cite{mccroryparusinski}
theorem 2.5.    
\end{proof}

\subsection{Arc-Symmetric Sets}
\label{AS} 

We call $X\subset \proj^N$ \emph{arc-symmetric} if 
for every real analytic arc $\gamma (t) : (-\varepsilon ,
\varepsilon) \to \proj^N$ such that $\gamma ((-\varepsilon ,0))
\subset X$  there is $\varepsilon ' >0$ such that 
$\gamma ((0,\varepsilon')) \subset X$.    
A subset $X\subset \R^N$ is called 
\emph{arc-symmetric} if it is arc-symmetric
as a subset of  $\proj^N$.  If $X$ is arc-symmetric 
then we write $X\in \AS$ for short.

\begin{rem}
The arc-symmetric sets were first introduced and studied by Kurdyka in
\cite{kurdyka1} but his original notion of arc-symmetric sets differs 
sligthly from ours.  In  \cite{kurdyka1} Kurdyka considers only 
closed semi-algebraic subsets 
$\R^N$ such that  
 for every real analytic arc $\gamma (t) : (-\varepsilon ,
\varepsilon) \to \R^N$ if $\gamma ((-\varepsilon ,0))
\subset X$ then 
$\gamma ((-\varepsilon ,\varepsilon)) \subset X$.  Note that this
definition does not give a sufficent control of the sets in 
question at infinty.  For instance one branch of a hyperbola is 
arc-symmetric in the sense of Kurdyka but not in ours.  
Its projective compactification is not arc-symmetric neither in
Kurdyka's sense nor in ours.  
But for the proof of surjectivity theorem, that is of
Corollary \ref{borelkurdyka}, Kurdyka uses in \cite{kurdyka2} 
the sets that remain 
arc-symmetric (in his sense) after compactification
of $\R^N$ that is in our terminology precisely the arc-symmetric
subsets of $\proj^N$ contained and closed in $\R^N$.  
Kurdyka's arc-symmetric sets are Euler, cf. \cite{mccroryparusinski},
but they do not form a constructible category. 
\end{rem}

The categories of algebraically constructible sets and 
arc-symmetric sets are very similar.  The topological properties 
of arc-symmetric sets can be studied by  means of Nash-constructible 
functions \cite{mccroryparusinski}, \cite{bonnard}, a tool analogous
 to algebraically constructible functions.  The following theorem 
shows these similarities.  It is simpler to us to formulate it 
for subsets of projective spaces so we avoid the passage to 
compactification.

\begin{thm}\label{ASequiv}
Let $X$ be a semi-algebraic subset of $\proj^N$.  
The following conditions are equivalent:  
\item [(i)]
$X$ is arc-symmetric.
\item [(ii)]
$\1 _X$ is Nash constructible as a function
(see \cite{mccroryparusinski}).
\item [(iii)] 
There exist blow Nash functions (see \cite{bonnard}) 
$g_1, \ldots, g_k$ such that $\1 _X = \sum \sgn g_i$.
\item [(iv)]
There exists a real algebraic morphism $F': Z'\to \proj^N$ and
$Z$ a union of connected components of $Z$ 
such that for $F= F'_{|Z}$ 
\begin{eqnarray*}
& x\in X \text { iff } \, \chi (F\inv (x)) = 1 \,\\
& x\notin X \text { iff } \, \chi (F\inv (x)) = 0 .
\end{eqnarray*} 
\item [(v)]
There exists a real algebraic morphism $F': Z'\to \proj^N$ and
$Z$ a union of connected components of $Z'$ 
such that for $F= F'_{|Z}$ 
\begin{eqnarray*}
& x\in X \text { iff } \, \chi (F\inv (x)) \equiv 1 \pmod 2\,\\
& x\notin X \text { iff } \, \chi (F\inv (x)) \equiv 0 \pmod 2 . 
\end{eqnarray*} 
\end{thm}

\begin{proof}
(i) is equivalent to (ii) by the argument of section 5 of 
\cite{mccroryparusinski}.
(ii) is equivalent to (iii) by
\cite{bonnard}.  (ii)  implies
(iv) by \cite{mccroryparusinski}.  (iv) implies (v)
trivially.

We show (v) implies (i).  For this we note that each closed
arc-symmetric subset of $\proj^N$ can be written uniquely
as a finite union of irreducible closed arc-symmetric sets,
see \cite{kurdyka1}.  On the other hand we have after
\cite{mccroryparusinski}.  

\begin{lem} 
Let $F':Z'\to \proj^N$ be a morphism of real algebraic
varieties, $Z\subset Z'$ closed and arc-symmetric.  Denote
$F= F'_{|Z'}$ and suppose that $F(Z) \subset Y$,
where $Y$ is closed arc-symmetric and irreducible.  Then there
exists a proper arc-symmetric subset $Y'\subset Y$, $\dim Y' <
\dim Y$,  such that 
$\chi (F\inv (x)) $ is constant modulo 2 on
$Y\setminus Y'$.   
\end{lem}

The rest of proof is similar to the algebraically constructible case.
\end{proof}

\begin{thm}\label{ASC}
Arc-symmetric  sets form a constructible category.
\end{thm}

\begin{proof} 
It follows from the topological properties of Nash-constructible
functions as  developed in  \cite{mccroryparusinski} section 5.  
First the links of Nash constructible functions are even that
 shows A4.  To show A3 we need to show the 
Nash constructible functions are stable by  inverse image 
and pushforward by mappings with arc-symmetric graph. 
It is clear from \cite{mccroryparusinski} section 5  the 
Nash constructible functions are stable by inverse image 
and pushforward by proper regular maps.  
Let $f:X\to Y$, $X\subset \proj^n$, $Y\subset \proj^m$, where 
$X,Y$ and the graph $\Gamma$ of $f$ are arc-symmetric.  
We denote the projections of $\proj^n\times \proj^m$ to 
$\proj^n$ and $\proj^m$ by 
$\pi_n$ and $\pi_m$ respectively.
Let $\varphi$, resp. $\psi$,  be a Nash constructible function 
supported in $X$, resp. $Y$.  On the level of constructible 
functions the restiction to $\Gamma$
correponds to the multiplication by $\1 _\Gamma$ and hence 
$f_*(\varphi) = (\pi_m)_* ((\pi_n)^* \cdot \1 _\Gamma )$, 
$f^*(\psi) = (\pi_n)_* ((\pi_m)^* \cdot \1 _\Gamma )$. 

 The rest of the proof is completely analogous to the proof of
 theorem \ref{ACC}. 
\end{proof}


\subsection{Remarks.}
\label{Remarks}
\medskip

Clearly $\AC$ is the smallest constructible category.
The difference between $\AC$ and $\AS$ lies in the fact
that a compact connected component of an arc-symmetric set is
arc-symmetric.  

\begin{prop}
Every constructible category 
is contained in $\AS$.  Furthermore, $\AS$ is the only 
constructible category containing the connected components
of compact real algebraic sets. 
\end{prop}

\begin{proof}
Let $\CC$ be a constructible category and let
$X\subset \proj^N$ be a $\CC$-constructible set.
We show that $X$ is arc-symmetric by induction on $\dim X$. 
Let $\dim X=n$ and suppose first that $X$ is closed in
$\proj^N$.   
Consider 
$$
\pi: \widetilde Y \to Y
$$
a resolution of $Y = \zar X$.  Let
$S\subset Y$ be an algebraic subset of $Y$ such that
$Y\setminus S$ is regular of pure dimension $n$,
$\dim S<n$, and
$$
\widetilde \pi := \pi_{|}: \widetilde Y \setminus
\pi \inv (S) \to Y\setminus S
$$
is an isomorphism.  Denote $\widetilde S := \pi \inv
(S)$.
$\pi \inv (X)$ is a $\CC$-constructible subset of 
$\widetilde Y$ and hence by lemma \ref{invariance}
\begin{equation}\label{union}
\pi \inv (X) = X_1 \cup X_2, 
\end{equation}
where $X_1$ is the union of some connected components of
$\widetilde Y$ and $\dim X_2<n$.
Note that $X' = \pi (X_1)
\cup S$ is arc-symmetric, see 
\cite{kurdyka1} for details. 
Let $S' = S\cup \zar{\pi (X_2)}$.  Then, by
the above $X\setminus S' = X'\setminus S'$ is arc symmetric.
On the other hand, $X\cap S'$ is $\CC$-constructible and
hence arc symmetric by the inductive assumption.  This
shows that $X$ is arc-symmetric.

Consider $\CC$-constructible $X\subset \proj^N$, not necessarily 
closed in $\proj^N$.  Then, by Remark \ref{remark}, 
$\dim (\cc X \setminus X) < \dim X$. $ \cc X$ is closed in 
$\proj^N$ and hence arc-symmetric by the above and 
$\cc X \setminus X$ is arc-symmetric by the inductive assumption. 
This shows that $X$ is arc-symmetric
and hence the first claim of proposition.

Suppose now that $\CC$ is a constructible category
containing the connected components
of compact real algebraic sets.  Let $X\subset \proj^N$ be closed 
and arc-symmetric.  We show that $X\in \CC$.
Apply to $X$ the construction of the first part of proof.
We find again \eqref{union}, where now $X_1\in \CC$.  Hence, 
be axiom A3, so does
$\pi (X_1\setminus \widetilde S)  = X \setminus S'$, where
again $S' = S\cup \zar{\pi (X_2)}$.  Again, $X\cap S'$ is
arc-symmetric and of dimension $< \dim X$ and hence belongs to
$\CC$ by the induction on dimension.  This shows that 
$X = (X\setminus S) \cup (X\cap S')$ is in $\CC$ as required.  
\end{proof}

By 
\cite{mccroryparusinski}  we have the following 
strenghtening of axiom 4. 

\begin{cor}
Every locally compact $\CC$-constructible set is Euler.
\end{cor}

\medskip
\section{Hierarchy of singularities}
\label{Hierarchy}

 We say that two germs $(Y,y)$  and $(X,x)$ of 
locally compact semi-algebraic sets are 
\emph{semi-algebraically homeomorphic}, and write 
$(Y,y)\sim (X,x)$,   
if there exists the germ of a semi-algebraic
homeomorphism $f:(Y,y) \to (X,x)$.   

Let $X$ be a locally compact semi-algebraic set.  
Then there is a semi-algebraic stratification of $X$ 
such that $(X,y)\sim (X,x)$ if 
$y,x$ belong to the same stratum.  Indeed, this follows from 
from the existence of a semi-algebraic triangulation of $X$.  
Moreover, by the same argument, this 
stratification can be chosen locally topologically trivial 
(i.e. for any stratum $S$ and any $x\in S$ there is 
a neighbourhood $U$ of $x$ in $X$ and a continuous retraction 
$\pi: U \to S\cap U$, $\pi|_{S\cap U} = id_{S\cap U}$,
such that $\pi$ is a trivial fibration). 
For each $x\in X$ consider 
$$
S_x = \{y\in X; (X,y)\sim (X,x)  \}.
$$
Clearly $S_x$ is a union of strata so it is semi-algebraic.  
Consequently at its generic point $S_x$ is a topological 
manifold.  But the germ of $X$ at arbitrary point of $S_x$ is 
homeomorphic to that at a generic point of $S_x$.  Thus we get 

\begin{prop}\label{strat}
Each $S_x$ is a locally closed semi-algebraic subset of 
$X$ and a pure-dimensional topological manifold.   
The decomposition into the equivalence classes 
of semi-algebraic homeomorphisms induces 
a canonical, locally topologically trivial, 
(topological) stratification of $X$.  
\end{prop}  

We call this stratification  \emph{the stratification by 
s.-a. topological types}.  

From now on we restrict ourselves to Euler locally compact 
semi-algebraic sets.  We denote by $\G$ the set 
of equivalence classes of such sets divided by semi-algebraic 
homeomorphisms.

\begin{defn}
Let $(Y,y)$ and $(X,x)$ be two germs of Euler locally compact 
 semi-algebraic sets.  We write $(Y,y)\prec (X,x)$ 
if there exists  the germ of  a continuous semi-algebraic 
injective map $i:(Y,y)\to (X,x)$.  
\end{defn}

\begin{prop}
If $(Y,y)\prec (X,x)$ and $(X,x)\prec (Y,y)$  then 
 $(Y,y)$ and $(X,x)$ are semi-algebraically homeomorphic.  
Thus $\prec$ gives a partial order on $\G$.  
\end{prop}

\begin{proof}  
$\prec$ is clearly transitive.  
If $(X,x)\prec (Y,y)$ and $(Y,y)\prec (X,x)$ then 
$(X,x)\prec (X,x)$ that is there exists a continuous 
semialgebraic injective map 
$i:( X, x)\to ( X,x)$ that must be a homeomorphism by 
proposition \ref{compact} (ii).  
\end{proof}

Note that after Lemma \ref{invariance} if $X$ is a topological 
manifold at $x$ of dimension $n$ and $(Y,y)\prec (X,x)$, 
with $\dim_y Y=n$, then $(Y,y)\sim (X,x)$.  
In a way $(Y,y)\prec (X,x)$ means that the singularity of 
$Y$ at $y$ is less complicated then that of $X$ at $x$ (at least 
if $\dim_y Y = \dim_x X$).  We want that the singularities 
around $x$ are not more complicated than that at $x$.

\begin{prop}\label{obok}
Let $X$ be an Euler locally compact  semi-algebraic set 
and let $x\in X$.  If $x\in \overline {S_y}$ and 
$(X,x)\prec (X,y)$ then $(X,x)\sim (X,y)$.  
\end{prop}

\begin{proof}
The proof is by induction on the codimension of $S_y$ in $X$.  
If this codimension is zero then $X$ is a topological manifold 
at $y$ and the result follows from corollary \ref{invariance2}.  

In the general case suppose that there exists a continuous
semi-algebraic injection $f:(X,x)\to (X,y)$.  Choosing $y'\in S_y$
arbitrarily close to $x$ we get $z = f(y')$ arbitrarily close to
$y$.  That is  
$y\in \overline {S_z}$ and $(X,y)\prec (X,z)$.  
Note that $y\in \overline {S_z}$ gives  $S_y \subset 
\overline {S_z}$, in particular, $\dim S_z > \dim S_y$.  
Therefore, by the inductive assumption, $(X,y)\sim (X,z)$.

\begin{lem}\label{along}
Let $X$ be an Euler locally compact semi-algebraic set, $S\subset X$ 
a connected component of a stratum of the stratification 
of $X$ by s.-a. topological types.  
Let $Y\subset X$ be Euler closed  semi-algebraic.  Then, if 
$$
(Y,x)=(X,x) 
$$
holds for an $x\in S$ then it does so for all $x\in S$.  
\end{lem} 

\begin{proof}
The set $S_1= \{x\in S; (Y,x)=(X,x)\}$ is clearly open in 
$S$.  We 
show that it is closed.     

Denote by $X_{reg}^d$ the set of points of $X$ at which 
$X$ is a topological manifold of dimension $d$.  Then 
$$
Y_{reg}^d \cap X_{reg}^d   
$$
is open and closed in $X_{reg}^d $.  Indeed it is clearly open.  
Suppose 
$x\in \overline {Y_{reg}^d \cap X_{reg}^d} \cap X_{reg}^d$.  
Then $x\in Y$ since $Y$ is closed in $X$ and $\dim_x Y$ must be $d$. 
Hence $x\in  Y_{reg}^d$ by lemma \ref{invariance}.  
This shows the claim.  

Let $x\in \overline S_1\cap S$.  Fix an integer $d$ and let $X'$ 
be a connected component of 
$X_{reg}^d$ such that $x\in \overline {X'}$.  
Hence, by local topological triviality of stratification, there 
is a neighbourhood $U_S$ of 
$x$ in $S$ such that $U_S\subset \overline {X'}$.  
$X'\cap Y_{reg}^d$ is nonempty since $U_S\cap S_1 \ne \emptyset$. 
Thus, by the above claim,  $X'\cap Y_{reg}^d = X'$ and hence 
$\overline {X'}\subset Y$.  Since the union of such $\overline
{X'}$ form a neighbourhood of $x$ in $X$ we get 
$(X,x)=(Y,x)$ as claimed.  
\end{proof}

We complete the proof of proposition \ref{obok}.  Let $B_y$ a small
open semi-algebraic neighbourhood (e.g. a ball) of $y$ in $X$ and
let $B_x = f\inv (B_y)$.  Then the restriction of $f$, 
$B_x\to B_y$ is injective and its image $f(B_x)$ is a closed
semi-algebraic and Euler subset of $B_y$.  
We apply lemma \ref{along} to $f(B_x)\subset B_y$ and the stratum
$S=S_y\cap B_y$.  Recall that there exists $z\in S$, $z=f(y')$, where
$y'\in S_y$ is close to $x$, such that $(X,y)\sim (X,z)$.  Then 
$f$ is a homeomorphism at $y'$ by proposition \ref{compact} (b). 
Thus $(f(B_x),z) = (B_y,z)$ and by lemma \ref {along} 
$(f(B_x),y) = (B_y,y)$.  That is $f$ is a local homeomorphism at 
$x$.  This ends the proof.  
\end{proof}

As a corollary we obtain the following

\begin{thm}\label{homeo1}
Let $X$ be a locally compact Euler semi-algebraic set of pure 
dimension $n$.  Let $f:X\to X$ be continuous and semi-algebraic. 
If $f$ is injective and surjective then it is a homeomorphism. 
\end{thm}

\begin{proof}
Let $\mathcal S$ be the stratification of $X$ by semi-algebraic 
topological types.  It suffices to show that for each stratum 
$S\in \mathcal S$ we have $f(S)=S$.  Indeed, then,  by proposition 
\ref{compact} (b), $f$ is a local homeomorphism at every $x\in X$
and the theorem follows.   The proof
is by induction on number of strata in $\mathcal S$. 
 If there is only 
one type, that is $X$ is a topological manifold, then the result 
follows from the invariance of domain or corollary \ref{invariance2}. 

Let $S$ be a stratum of $\mathcal S$ of maximal singularity type 
(with respect to $\prec$).   Then $f(S) \subset S$.  By 
lemma \ref{invariance}, $f(S)$ is 
open in $S$.
We show that $f(S)$ is also closed in $S$.  
By propositions \ref{strat} and \ref{obok},  
$S$ is a topological manifold closed in $X$.  Consider the restriction 
 of $f$, $f\inv (S) \to S$.  
Let $y\in S$, $f(x) =y$,  and choose, as in the proof of 
proposition \ref{obok},  $B_y$ a small
open semi-algebraic neighbourhood of $y$ in $X$ and let 
$B_x = f\inv (B_y)$.  Then, see the proof of lemma \ref{along}, 
$Y' = (f(B_x))_{reg} \cap (B_y)_{reg}$ is open and closed in 
$(B_y)_{reg}$.  We have $y\in \overline {Y'}$ and hence, 
by local topological triviality of stratification $\mathcal S$, 
$S\cap B_y 
\subset \overline {Y'}$.  This shows that the mapping 
$f\inv (S) \to S$ is open. 
Consequently, since $S$ is closed in $X$, 
$S\setminus f(S) = f(f\inv (S)\cap (X\setminus S))$ is open in
$S$.  Thus we have showed that $f(S)$ is open and closed in $S$ 
that is a union of connected components of $S$.  
Corollary \ref{invariance2} shows that the restriction of
$f$, $S\to f(S)$ is a homeomorphism, in particular $S$ and $f(S)$
are homeomorphic.   Hence $f(S) =S$.   

To complete the proof of theorem we apply the inductive
assumption to $f|_{X\setminus S} : X\setminus S \to X\setminus S$.  
\end{proof}

\begin{thm}\label{homeo2}
Let $X$ be a locally compact $\CC$-constructible 
semi-algebraic set and let $f:X\to X$ be an injective continuous 
$\CC$-constructible map.  Then $f$ is a homeomorphism.  
\end{thm}

\begin{proof}
$f$ is surjective by theorem \ref{mainthm}.  
Let $\dim X=n$ and denote by $\mathcal S$ be
 the stratification of $X$ by semi-algebraic 
topological types.  As in the proof of theorem \ref{homeo1} it suffices
to show that $f(S)=S$ for each stratum $S$.  

\begin{lem}
There is a filtration of $X$ by
 closed in $X$, $\CC$-constructible, and stable by $f$ sets 
$$
X=X^n \supset X^{n-1}\supset \cdots \supset X^0 \supset X^{-1}= 
\emptyset
$$
such that for every $i$:  $\cir X^i:=X^i\setminus X^{i-1}$ is a 
topological 
manifold of pure dimension $i$ (or empty) and 
$(X,y)\sim (X,x)$ for all $y,x$
belonging to the same connected component of $\cir X^i$.  

(stability by $f$ for $X^i$ means $f(X^i)\subset X^i$. By theorem
\ref{mainthm} it is equivalent to $f(X^i)= X^i$.)
\end{lem}

\begin{proof}
By descending induction on $k$ we construct a filtration 
\begin{equation}\label{shortfiltration}
X=X^n \supset X^{n-1}\supset \cdots \supset X^k 
\end{equation}
satisfying all the above properties except that $X^k$ can be 
arbitrary closed in $X$, $\CC$-constructible, and stable by $f$ 
set of dimension at most $k$.  If $k=n-1$ then we 
take as $X^{n-1}$ the 
$\CC$-closure of the singularities of $X$, i. e. the set 
 $\tilde \Sigma$ of the proof of theorem \ref{borelkurdyka}.  
Suppose that the filtration \eqref{shortfiltration} has been already 
constructed.  We construct $X^{k-1}$.  Set 
$$
Y = \bigcup_{S\in\mathcal S} \Sigma (S\cap X^k) 
$$
where $\Sigma (S\cap X^k)$ denote the set of points where  
$S\cap X^k$ is not topological manifold of
dimension $k$.   Then $X^k\setminus Y$ is a topological manifold
of dimension k.  Let  
$$
X^{k-1} = \cc {Y}.
$$
It is clear that $X^{k-1}$ satisfies the required properties.  
\end{proof}

Fix an integer $i$.  
By construction $f$ sends the connected components of $\cir X^i$
onto the connected components of $\cir X^i$.  To each such 
connected component corresponds one singularity type.  Thus the 
components corresponding to the maximal types (with respect to 
$\prec$) are sent onto the components of the same type.  
Hence, by inductive argument on the number of types (or components),
we see that for each stratum $S$, $f(S\cap \cir X^i) = (S\cap \cir
X^i)$.  
This, for all $i$, shows $f(S) =S$.  This ends the proof of
theorem.  
\end{proof}

\begin{cor}\label{borelkurdyka2} 
Any injective regular 
self-mapping of a real algebraic variety is a homeomorphism.
\end{cor}

\begin{cor} 
Any injective continuous self-mapping of a locally 
compact arc-symmetric set, 
the graph of which is arc-symmetric, is a homeomorphism.
\end{cor}


\end{document}